\RequirePackage[english]{babel}
\documentclass[a4paper,twoside,11pt]{amsart}
\usepackage[latin1]{inputenc}
\usepackage{amsmath}
\usepackage{amssymb}
\usepackage{amsfonts}
\usepackage{mathrsfs}
\usepackage{fixltx2e}[2005/12/01]
\usepackage{hyperref}

\usepackage{enumerate}

\usepackage[matrix,arrow,cmtip,frame,graph,2cell]{xy} 
\SelectTips{cm}{}
\newdir{(}{{}*!/-5pt/@^{(}}
\newdir{(x}{{}*!/-5pt/@_{(}}
\newdir{+}{{}*!/-9pt/{}}
\newdir{>+}{@{>}*!/-9pt/{}}
\entrymodifiers={+!!<0pt,\fontdimen22\textfont2>}
\UseAllTwocells

\usepackage[nosubsections]{mytheorems2}

\newcommand\pref[1]{\eqref{#1}}

\newcommand\map[3]{#1\colon #2\rightarrow #3}
\newcommand\inj{\hookrightarrow}
\newcommand\iso{\cong} 
\newcommand\id[1]{\mathrm{id}_{#1}} 
\DeclareMathOperator{\Hom}{Hom}
\newcommand\sO{\mathcal{O}}
\DeclareMathOperator{\Spec}{Spec}
\newcommand{\etale}{\'{e}tale}
\newcommand\SG[1]{{\mathfrak{S}_{#1}}}   
\newcommand{\catHom}{\mathbf{Hom}}

\newcommand\unram{\mathrm{unram}}
\newcommand\qfin{\mathrm{qfin}}
\newcommand\repr{\mathrm{repr}}

\newcommand{\stU}{\mathscr{U}}
\newcommand{\stX}{\mathscr{X}}
\newcommand{\stY}{\mathscr{Y}}
\newcommand{\stZ}{\mathscr{Z}}
\newcommand{\stHom}{\mathscr{H}om}
\newcommand\FHilb{\mathcal{H}\mathit{ilb}}
\newcommand\Hilb{\mathrm{Hilb}}
\newcommand{\AF}{{\upshape AF}}
\newcommand\HilbSt{\mathscr{H}}
\newcommand\EtSt{\mathscr{E}t}
\newcommand\Et{\mathrm{ET}}
\newcommand\Sec{\mathrm{SEC}}
\newcommand\weilr{\mathbf{R}} 

\newcommand\BG{\mathrm{B}G}
\newcommand\BGL{\mathrm{BGL}}

\newcommand\Set{\mathbf{Set}}
\newcommand\Sch{\mathbf{Sch}} 

\newcommand{\equalizer}[2]{\xymatrix@1@M=0mm@C=7mm{#1%
\ar@<.5ex>@{+->+}[r] \ar@<-.5ex>@{+->+}[r] & #2}}

\begin{document}

\title[Hilbert schemes and Hilbert stacks of points]
{Representability of Hilbert schemes and Hilbert stacks of points}
\author{David Rydh}
\address{Department of Mathematics, University of California, Berkeley,
970 Evans Hall \#3840, Berkeley, CA 94720-3840 USA}
\email{dary@math.berkeley.edu}
\date{2011-07-09}
\thanks{Supported by the Swedish Research Council 2008-7143.}


\begin{abstract}
  We show that the Hilbert functor of points on an arbitrary \emph{separated}
  algebraic stack is an algebraic space. We also show the algebraicity of the
  Hilbert stack of points on an algebraic stack and the algebraicity of the
  Weil restriction of an algebraic stack along a finite flat morphism.
  For the latter two results, no separation assumptions are necessary.
\end{abstract}

\maketitle


\setcounter{secnumdepth}{0}
\begin{section}{Introduction}
  The purpose of this note is to give a short and elementary proof of the
  algebraicity of the Hilbert functor of points $\FHilb^d_{X}$ for an arbitrary
  separated algebraic space or stack~$X$,
  cf.\ Theorem~\pref{T:Hilb-functor-alg}.
  The key fact is that given an \etale{} (resp.\ smooth) presentation
  $\map{f}{U}{X}$ there is an open subspace $\FHilb^d_{U\to X}$ of
  $\FHilb^d_{U}$ and a representable \etale{} (resp.\ smooth) surjective map
  $\map{f_*}{\FHilb^d_{U\to X}}{\FHilb^d_X}$.
  It follows that $\FHilb^d_{X}$ is algebraic if $\FHilb^d_{U}$ is algebraic
  and for affine $U$ this is
  well-known~\cite{nori_Hilb_appendix,GLS_Affine_Hilb}.

  When $X$ is not separated, the Hilbert functor is not
  representable~\cite{lundkvist-skjelnes_noneff-deform}. A replacement for the
  Hilbert functor is then the Hilbert stack which was briefly introduced by
  M.\ Artin~\cite[App.]{artin_versal_def_alg_stacks}. Applying the same
  method as for the Hilbert functor, we deduce the algebraicity of the Hilbert
  stack of points from the affine case,
  cf.\ Theorems~\pref{T:alg-of-HilbSt-aff}
  and~\pref{T:algebraicity-of-HilbSt}.
  Along the way we also prove the algebraicity of the Weil restriction of an
  algebraic stack along a finite flat morphism and the algebraicity of
  Hom-stacks when the source is finite flat, cf.\ Theorems~\pref{T:alg-of-weilr}
  and~\pref{T:alg-of-Hom}.
  In Section~\ref{S:etale-families}, we show that the open substack
  parameterizing \etale{} families coincides with the stack quotient of the
  $d$\textsuperscript{th} fiber product $X^d$ by the symmetric group. This is
  the stack of zero-dimensional
  branchvarieties~\cite{alexeev-knutson_branch-varieties}.

  Let us indicate the relationship between our methods and more standard
  representation techniques. Thus consider the following diagram of
  subcategories of algebraic spaces over some fixed base scheme:
\vspace{-3mm}
{
\newcommand{\cQProj}{\textbf{QProj}}
\newcommand{\cAff}{\textbf{Aff}}
\newcommand{\fp}{\textrm{fp}}
\newcommand{\sep}{\textrm{sep}}
\newcommand{\cAF}{\textbf{AF-Sch}} 
\newcommand{\cAlgSp}{\textbf{AlgSp}} 
$$\xygraph{ !{0;<3.5cm,0cm>:<0cm,1.5cm>::}
[d]
*++{\cAff}="AFF" :@{(->} [r]
*++{\cAF}="AF" [u]
*++{\cQProj_{\fp}}="QP"
 :@{(->} "AF" "QP" :@{(->} [r]
*++{\cAlgSp_{\fp,\sep}}="SASFP" :@{(->} [d]
 *++{\cAlgSp_{\sep}}="SAS"
"AF" :@{(->} "SAS" "SASFP" :@{(->} [r]
*++{\cAlgSp_{\fp}} :@{(->} [d]
 *++{\cAlgSp.}="AS"
"SAS" :@{(->} "AS"
}$$
Here $\cAff$ (resp.\ $\cQProj$, resp.\ $\cAF$, resp.\ $\cAlgSp$) denotes the
category of affine schemes (resp.\ quasi-projective schemes, resp.\
\AF{}-schemes, resp.\ algebraic spaces). A scheme is \AF{} if every
finite set of points lies in an affine open subscheme. The subscripts ``$\fp$''
and ``$\sep$'' stand for ``locally of finite presentation'' and ``separated''.

Let $X$ be an algebraic space in one of these categories and consider a functor
or stack of objects on $X$ such as the Hilbert functor $\FHilb_{X}$. The
classical approach is to work in the category $\cQProj_\fp$ and show that the
functor can be embedded into a projective scheme~\cite[No.~221]{fga}.
Alternatively, we can work in the bigger category $\cAlgSp_{\fp,\sep}$ where
Artin's algebraization
theorems~\cite{artin_alg_formal_moduli_I,artin_versal_def_alg_stacks} apply.
In this category, the algebraicity of the Hilbert functor, Hilbert stack and
Weil restriction is
well-known, cf.\ \cite[Cor.~6.2]{artin_alg_formal_moduli_I},
\cite[Rmk.~4.5]{starr_artins_axioms} and~\cite[\S2.1]{lieblich_Coh_stack}.

In this paper we will primarily be interested in the second line of the diagram
and show that in the \emph{zero-dimensional case} we can use more elementary
and constructive methods to deduce algebraicity. As a bonus, we need neither
finiteness nor separatedness assumptions. Similar methods were applied
in~\cite{GLS_Affine_Hilb,GLS_Affine_Quot} for Hilbert and Quot schemes of
affine and \AF{}-schemes.

In a subsequent paper, we will turn the attention to the category
$\cAlgSp_{\fp}$ and use Artin's algebraization theorem
to show the algebraicity of the Hilbert stack parameterizing
\emph{higher-dimensional} families on a given, possibly \emph{non-separated},
scheme or algebraic space. This relies on variants of Chow's lemma and
Grothendieck's existence theorem for non-separated spaces.
}

As it is easier to first establish the algebraicity of
the Hilbert stack and then \emph{a posteriori} verify that it is
quasi-separated we work with general algebraic spaces and algebraic stacks
without any separation assumptions, cf.\ Appendix~\ref{A:stacks}.

\begin{subsection}{Acknowledgments}
I would like to thank D.\ Laksov, R.\ Skjelnes and the referee for useful
comments and discussions.
\end{subsection}

\end{section}
\setcounter{secnumdepth}{3}


\begin{section}{The Hilbert functor and the Hilbert stack}
For simplicity, we work over a fixed base \emph{scheme} $S$. If $X$ is
a scheme or an algebraic stack over $S$, then a property of $X/S$ always refers
to a property of the structure morphism $X\to S$.

\begin{definition}
We say that a finite morphism $\map{f}{X}{Y}$ is \emph{flat of rank~$d$} if
$f_*\sO_X$ is a locally free $\sO_Y$-module of constant rank~$d$.
\end{definition}

\begin{definition}
Let $X/S$ be a \emph{separated} scheme (resp.\ separated algebraic space,
resp.\ separated algebraic stack). The \emph{Hilbert functor of points}
$\FHilb^d_{X/S}$ is the functor which to an $S$-scheme $T$ assigns the set of
closed subschemes (resp.\ subspaces, resp.\ substacks) $Z\inj X\times_S T$
such that the second projection $\map{p}{Z}{T}$ is finite and flat of
rank~$d$.
\end{definition}

\begin{remark}
It is easily seen, using~\cite[Thm.~12.2.1\ (i),\ (ii)]{egaIV}, that
$\FHilb^d_{X/S}$ is an open and closed subfunctor of the full Hilbert functor
$\FHilb_{X/S}$ which parameterizes closed subspaces (or substacks) which are
flat, proper and of finite presentation. Note that an object of
$\FHilb^d_{X/S}(T)$ is a closed substack $Z\inj X\times_S T$ such that
$\map{p}{Z}{T}$ is finite and hence $Z$ is always a \emph{scheme} even if $X$
is an algebraic stack.
\end{remark}

\begin{remark}
Using Artin's criteria for algebraicity it can be shown that the Hilbert
functor $\FHilb_{\stX/S}$ is a separated algebraic space locally of finite
presentation when $\stX/S$ is a separated algebraic stack locally of finite
presentation. In this generality, the result is due to
M.\ Olsson~\cite[Thm.~1.5]{olsson_proper-coverings} but also
see~\cite{olsson-starr_quot-functors} for boundedness results when $\stX$ is a
Deligne--Mumford stack.
\end{remark}

We will now define the Hilbert \emph{stack} of
points~\cite[App.]{artin_versal_def_alg_stacks}. The difference between
the Hilbert stack and the Hilbert functor is that in the stack we consider
flat families $Z\to T$ with morphisms $Z\to X$ without the condition that
$Z\to X\times_S T$ is a closed immersion.

\begin{definition}
Given an algebraic stack $\stX/S$, let $\HilbSt^d_{\stX}$ be the category
with objects pairs of morphisms $({\map{p}{Z}{T}},\map{q}{Z}{\stX})$
where $T$ is an $S$-scheme and $p$ is finite and flat of rank~$d$. The
morphisms are triples $(\varphi,\psi,\tau)$ fitting into a $2$-commutative
diagram
\vspace{-7 mm}
$$\xymatrix{
{Z_1}\ar[r]^{\varphi}\ar[d]^{p_1}\rruppertwocell<9>^{q_1}{<-2.8>\tau}
  & {Z_2}\ar[d]^{p_2}\ar[r]^{q_2} & \stX \\
{T_1}\ar[r]^{\psi} & {T_2}\ar@{}[ul]|\square
}$$
such that the square is cartesian. The category $\HilbSt^d_{\stX}$ is fibered in
groupoids over $\Sch_{/S}$ and by \etale{} descent of affine
schemes~\cite[Exp.~VIII, Thm.~2.1]{sga1} it follows that $\HilbSt^d_{\stX}$ is a
stack. We call $\HilbSt^d_{\stX}$ the \emph{Hilbert stack of $d$ points} on~$\stX$.
\end{definition}

When $\stX/S$ is a \emph{separated} algebraic stack, then the Hilbert functor
$\FHilb^d_{\stX/S}$ is an open substack of the Hilbert stack $\HilbSt^d_{\stX}$,
cf.\ Proposition~\pref{P:substacks-of-Hilb} below.

\begin{example}
The Hilbert stack $\HilbSt^1_{\stX}$ is equivalent to $\stX$. If $\stX$ is
separated, then the Hilbert functor $\FHilb^1_{\stX}$ equals the
automorphism-free locus of $\stX$.
\end{example}

In the remainder of this section we will review more general Hilbert stacks
present in the literature. These are not used in the subsequent sections.

\begin{definition}
Let $\stX$ be an algebraic stack. The Hilbert stack $\HilbSt_{\stX}$ is the
stack parameterizing flat and proper algebraic stacks $\map{p}{\stZ}{T}$ of
finite presentation together with a representable morphism
$\map{q}{\stZ}{\stX}$.
%
%
The substack of objects such that $\map{(q,p)}{\stZ}{\stX\times_S T}$ is
locally quasi-finite (resp.\ unramified) is denoted
$\HilbSt^{\qfin}_{\stX}$ (resp.\ $\HilbSt^{\unram}_{\stX}$). The substack of
objects such that $p$ is representable, i.e., $\stZ$ is an algebraic space, is
denoted $\HilbSt^{\repr}_{\stX}$.
\end{definition}

\begin{lemma}\label{L:locus-lemma}
Let $\map{f}{X}{Y}$ be a morphism of algebraic spaces, locally of finite type.
\begin{enumerate}
\item The locus of points in $X$ where $f$ is quasi-finite is open.
\item The locus of points in $X$ where $f$ is unramified is open.
\item If $f$ is proper, then the subfunctor $Y'\subseteq Y$ consisting of
morphisms $T\to Y$ such that $X\times_Y T\to T$ is a closed immersion, is
an open subspace.\label{LI:closed-imm}
\end{enumerate}
\end{lemma}
\begin{proof}
(i) is~\cite[Cor.~13.1.4]{egaIV} as the question is \etale{}-local on $X$ and
$Y$. The locus where $f$ is unramified is the complement of the support
of $\Omega^1_f$ and thus open.
To show (iii) we can assume that $f$ is quasi-finite and hence
finite~\cite[Cor.~A.2.1]{laumon} since the locus of $Y$ where $f$ is
quasi-finite is open by (i). That $Y'\subseteq Y$ is an open subspace
now follows by Nakayama's lemma.
\end{proof}

\begin{proposition}\label{P:substacks-of-Hilb}
Let $\stX/S$ be an algebraic stack. Then
\begin{enumerate}
\item $\HilbSt^{\repr}_{\stX}\subset \HilbSt_{\stX}$ is an open substack.
\item $\HilbSt^d_{\stX}\subset \HilbSt^{\repr}_{\stX}$ is an open and closed
substack.
\item $\HilbSt^{\unram}_{\stX}\subset \HilbSt^{\qfin}_{\stX}\subset \HilbSt_{\stX}$ are open substacks.
\item If $\stX/S$ is a \emph{separated} algebraic stack, then the Hilbert
functor $\FHilb_{\stX}$ is an open subfunctor of the Hilbert stack
$\HilbSt_{\stX}$.
\end{enumerate}
\end{proposition}
\begin{proof}
Let $T$ be a scheme and let $(\map{p}{\stZ}{T},\map{q}{\stZ}{\stX})$ be an
object of $\HilbSt_{\stX}(T)$. We let
$T^{\repr}=T\times_{\HilbSt_{\stX}}\HilbSt^{\repr}_{\stX}\subseteq T$ and
similarly for $T^{\unram}$ and $T^{\qfin}$.

Note that since $p$ is separated, the inertia stack $I_\stZ\to \stZ$ is
proper. By Lemma~\pref{L:locus-lemma}~\ref{LI:closed-imm} applied to
$I_\stZ\to \stZ$ the automorphism-free locus $\stZ'\subseteq |\stZ|$ is
open. Thus $T^{\repr}=T\setminus p(\stZ\setminus \stZ')$ is an open subscheme
of $T$. This shows (i). The second statement follows
from~\cite[Thm.~12.1.1\ (i),\ (ii)]{egaIV}.

To show (iii), let $\stZ'\subseteq \stZ$ be the locus where
$\map{(q,p)}{\stZ}{\stX\times_S T}$ is locally quasi-finite
(resp.\ unramified). This locus is open by Lemma~\pref{L:locus-lemma} and thus
$T^{\qfin}$ (resp.\ $T^{\unram}$) is the open subscheme $T\setminus
p(\stZ\setminus \stZ')$.

(iv) If $\stX/S$ is separated, then $\map{(q,p)}{\stZ}{\stX\times_S T}$ is
proper and by Lemma~\pref{L:locus-lemma}~\ref{LI:closed-imm}, the locus of
$\stX\times_S T$ over which $(q,p)$ is a closed immersion is open. We let
$\stZ'\subseteq \stZ$ be the inverse image of this locus. Then
$T\times_{\HilbSt_{\stX}}\FHilb_{\stX}=T\setminus p(\stZ\setminus \stZ')$ is an
open subscheme of $T$.
\end{proof}

The stack $\HilbSt_{\stX}$ is not algebraic in general as Grothendieck's
existence theorem does not hold without any projectivity assumptions on either
$p$ or $(q,p)$.
The algebraicity of $\HilbSt^{\qfin,\repr}_{\stX}$ for a \emph{separated}
algebraic stack $\stX$ with finite diagonal, locally of finite presentation
over $S$, is proved by J.\ Starr in~\cite[Rmk.~4.5]{starr_artins_axioms}.
A.\ Vistoli studies the substack $\HilbSt^{\unram,\repr}_{\stX}$ when $\stX$
is a separated Deligne--Mumford stack, locally of finite presentation over $S$
although he does not prove algebraicity~\cite{vistoli_hilbert_stack}.

A sketch of the proof of the algebraicity of $\HilbSt^{\qfin}_{\stX}$ for a
\emph{separated} algebraic stack, locally of finite presentation over $S$,
is given by M.\ Lieblich in~\cite[Thm.~2.1, Cor.~2.5]{lieblich_Coh_stack}.
More generally, it is expected that $\HilbSt^{\qfin}_{\stX}$ is algebraic when
$\stX/S$ is an algebraic stack locally of finite presentation with separated
and quasi-finite diagonal. Note that in this case
$\map{(q,p)}{\stZ}{\stX\times_S T}$ is quasi-finite and separated, hence
quasi-affine so that $\sO_\stZ$ is $(q,p)$-ample.
%
%


The \emph{stack of branchvarieties}, introduced by V.\ Alexeev and
A.\ Knutson~\cite{alexeev-knutson_branch-varieties}, is the open substack of
$\HilbSt^{\qfin}_{\stX}$ parameterizing families $\stZ\to T$ with geometrically
reduced fibers.

The \emph{space of husks} of a separated algebraic space $X$ is the open
subspace of $\HilbSt^{\qfin}_{X}$ parameterizing families $Z\to T$ with
$\left(S_1\right)$-fibers such that $\map{(q,p)}{Z}{X\times_S T}$ is
birational onto its image. This is the space of algebra husks of quotients of
$\sO_X$ introduced by
J.\ Koll\'ar~\cite[Rmk.~10]{kollar_simultaneous-norm-and-alg-husks}. When
limited to $1$-dimensional husks this is M.\ H\o{}nsen's space of
Cohen--Macaulay curves~\cite{honsen_mit-thesis}.

\end{section}


\begin{section}{Representability of the Hilbert scheme}
In this section we show that the Hilbert functor $\FHilb^d_{X}$ is separated
and that it is represented by a scheme if $X$ is an \AF{}-scheme. As we do not
know \emph{a priori} that $\FHilb^d_{X}$ is quasi-separated, it is not enough
to check the valuative criterion for separatedness.

\begin{lemma}\label{L:Hilb-separated}
Let $S$ be a scheme and let $\stX\to S$ be a separated algebraic stack. Then
the diagonal of $\FHilb^d_{\stX/S}\to S$ is a closed immersion.
\end{lemma}
\begin{proof}
Let $T$ be an $S$-scheme and let $T\to \FHilb^d_{\stX/S}\times_S
\FHilb^d_{\stX/S}$ be a morphism corresponding to closed subschemes $Z_1\inj
\stX\times_S T$ and $Z_2\inj \stX\times_S T$.  Let
$\map{T_{12}}{\Sch_{/T}}{\Set}$ be the functor defined as follows: for a
$T$-scheme $T'$, we let $T_{12}(T')$ be the one-point set if $Z_1\times_{T}
T'=Z_2\times_{T} T'$ and the empty set otherwise. Then $T_{12}\inj T$ is
represented by a closed subscheme~\cite[Lem.~9.7.9.1]{egaI_NE}. Since the
natural diagram
$$\xymatrix@C+5mm{\FHilb^d_{\stX/S}\ar@{(->}[r]^-{\Delta}
 & \FHilb^d_{\stX/S}\times_S \FHilb^d_{\stX/S}\\
T_{12}\ar@{(->}[r]\ar[u] & T\ar[u]_{[Z_1],[Z_2]}}$$
is cartesian, it follows that $\Delta_{\FHilb^d_{\stX/S}}$ is a closed
immersion.
\end{proof}

\begin{definition}
Let $X$ be a scheme. We say that $X$ is an \AF{}-scheme if every finite set
of points $Z\subseteq X$ is contained in an affine open subset of~$X$.
We say that a morphism of stacks $\map{f}{\stX}{\stY}$ is \AF{} if for every
affine scheme $T$ and morphism $T\to \stY$, the fiber product
$\stX\times_{\stY} T$ is an \AF{}-scheme.
\end{definition}

\begin{remark}\label{R:lqf+sep-AF}
If $S$ is an affine scheme and $X\to S$ is a locally quasi-finite and separated
morphism of algebraic spaces, then it follows from Zariski's main
theorem~\cite[Thm.~A.2]{laumon} that every finite subset $Z\subseteq X$ is 
contained in a quasi-affine open subscheme of $X$. It then follows
from~\cite[Cor.~4.5.4]{egaII} that $X$ is an \AF{}-scheme.
In particular, every representable, locally quasi-finite and separated morphism
is \AF{}.
\end{remark}

\begin{theorem}\label{T:repr-of-Hilb-AF}
Let $S$ be an affine scheme and let $X/S$ be an \AF{}-scheme. Then
$\FHilb^d_{X/S}$ is represented by a separated scheme $\Hilb^d(X/S)$.
\end{theorem}
\begin{proof}
When $X$ is an affine scheme then it is known that $\FHilb^d_{X/S}$ is
represented by a scheme~\cite{nori_Hilb_appendix,GLS_Affine_Hilb}.
If $X$ is an \AF{}-scheme, then let $X=\bigcup_\alpha U_\alpha$ be an open
cover of $X$ by affines such that every subset of $d$ points of $X$ lies in
some $U_\alpha$. It is then easily seen that $\coprod_\alpha
\FHilb^d_{U_\alpha}\to \FHilb^d_X$ is a Zariski covering and thus
$\FHilb^d_X$ is represented by a scheme.
%
%
That $\Hilb^d(X/S)$ is separated is Lemma~\pref{L:Hilb-separated}.
\end{proof}

\end{section}


\begin{section}{Weil restriction and push-forward of Hilbert stacks}
In this section we show that the \emph{Weil restriction} of an algebraic space
or algebraic stack along a finite flat morphism is algebraic. As a consequence
the \emph{push-forward} of Hilbert stacks is algebraic. We also provide a
list of properties for the Weil restriction and the push-forward of Hilbert
stacks. This section generalizes the results of Bosch, L\"utkebohmert and
Raynaud~\cite[\S7.6]{raynaud-bosch-lutkebohmert_Neron-models}.

\begin{definition}
Let $X\to S'$ and $\map{g}{S'}{S}$ be morphisms of algebraic spaces. The
\emph{Weil restriction} $\weilr_{S'/S}(X)$ is the functor from
$S$-schemes to sets that takes an $S$-scheme $T$ to the set of sections of
$X_T\to S'_T$, i.e.,
$$\weilr_{S'/S}(X)(T)=\Hom_{S'}(S'\times_S T, X).$$
Similarly, if $\stX\to S'$ is an algebraic stack, then
$\weilr_{S'/S}(\stX)$ is the stack with $T$-points the groupoid
$\catHom_{S'}(S'\times_S T, \stX)$
and the natural notion of pull-back. If $\map{f}{\stX}{\stY}$ is a morphism
of algebraic stacks, then there is a natural morphism of stacks
$$\map{\weilr_{S'/S}(f)}{\weilr_{S'/S}(\stX)}{\weilr_{S'/S}(\stY)}$$
taking a morphism $\map{s}{S'\times_S T}{\stX}$
to the morphism $\map{f\circ s}{S'\times_S T}{\stY}$.
\end{definition}


The Weil restriction is also sometimes denoted $\Pi_{S'/S} X$ or $g_*X$ and
is also known as \emph{restriction of scalars}, cf.~\cite[No.~195, \S C 2]{fga}
and~\cite{olsson_Hom-stacks}. The Weil restriction can also be defined on
$2$-morphisms so that $\weilr_{S'/S}$ becomes a strict $2$-functor from the
$2$-category of stacks over $S'$ to the $2$-category of stacks over $S$.
This functor is \emph{left exact}, i.e., takes $2$-fiber products to
$2$-fiber products and the terminal object $S'$ to the terminal object $S$.
Indeed, $\weilr_{S'/S}$ is a right $2$-adjoint to the pull-back functor
$g^{-1}$.

\begin{definition}
Let $\map{f}{\stX}{\stY}$ be a morphism of algebraic stacks. There is a
natural morphism $\map{f_*}{\HilbSt^d_{\stX}}{\HilbSt^d_{\stY}}$ taking an
object $(\map{p}{Z}{T},\map{q}{Z}{\stX})$ to $(p,f\circ q)$.
\end{definition}

\begin{definition}
Let $P$ be a property of morphisms of stacks. We say that $P$ is \emph{stable
under base change} if for every morphism $\map{f}{\stX}{\stY}$ with $P$ and
every morphism $\stY'\to \stY$ the base change $\map{f'}{\stX\times_{\stY}
\stY'}{\stY'}$ has $P$. We say that $P$ \emph{can be checked on affines}, if
a morphism $\map{f}{\stX}{\stY}$ has $P$ if and only if
$\map{f}{\stX\times_{\stY} T}{T}$ has $P$ for every affine scheme $T$ and
morphism $T\to \stY$.
\end{definition}

A stable property which is fppf-local on the base can be checked on affines. An
example of a non-fppf local property that can be checked on affines is the
property ``strongly representable'', i.e., represented by schemes.
Another example is the property \AF{}.

\begin{lemma}\label{L:property-equivalence}
Let $P$ and $Q$ be two properties of morphisms of stacks. Assume that these
properties are stable under base change and that $Q$ is a property that can be
checked on affines. Let $d$ be a positive integer. The following are equivalent:
\begin{enumerate}
\item If $S'\to S$ is a finite flat morphism of rank~$d$ between
algebraic spaces and if $\map{f}{\stX}{\stY}$ is a morphism of algebraic
stacks over $S'$ with property $P$ then
$\map{\weilr_{S'/S}(f)}{\weilr_{S'/S}(\stX)}{\weilr_{S'/S}(\stY)}$ has
property $Q$.\label{LI:pe:weilr-map}
\item If $S$ is an affine scheme, if $S'\to S$ is finite flat of rank
$d$ and if $\stX\to S'$ is a morphism of stacks with property $P$, then
$\weilr_{S'/S}(\stX)\to S$ has property $Q$.\label{LI:pe:weilr}
\item If $\map{f}{\stX}{\stY}$ is a morphism of algebraic stacks over $S$ with
property $P$ then $\map{f_*}{\HilbSt^d_{\stX/S}}{\HilbSt^d_{\stY/S}}$ has
property $Q$.\label{LI:pe:pushfwd}
\end{enumerate}
\end{lemma}
\begin{proof}
We have that \ref{LI:pe:weilr-map}$\implies$\ref{LI:pe:weilr} since
$S=\weilr_{S'/S}(S')$.

\ref{LI:pe:weilr}$\implies$\ref{LI:pe:pushfwd}:
Let $\map{f}{\stX}{\stY}$ be a morphism of algebraic $S$-stacks, let $T$ be an
affine $S$-scheme and let $\map{h}{T}{\HilbSt^d_{\stY/S}}$ be a morphism
corresponding to a family $(Z\to T,Z\to \stY)$. Then the diagram
\begin{equation}\label{E:pushfwd-weil-cart}
\vcenter{\xymatrix{%
\weilr_{Z/T}(\stX\times_{\stY} Z)\ar[r]\ar[d] & T\ar[d]^{h} \\
\HilbSt^d_{\stX/S}\ar[r]^{f_*} & \HilbSt^d_{\stY/S}\ar@{}[ul]|\square}}
\end{equation}
is cartesian and thus~\ref{LI:pe:weilr}$\implies$\ref{LI:pe:pushfwd}.

\ref{LI:pe:pushfwd}$\implies$\ref{LI:pe:weilr-map}:
Let $\map{g}{S'}{S}$ be finite flat of rank~$d$ and
let $\map{f}{\stX}{\stY}$ be a morphism of algebraic $S'$-stacks.
Let $\map{h}{S}{\HilbSt^d_{S'/S}}$ be the morphism corresponding to the family
$(g,\id{S'})$. Then from the cartesian diagram~\eqref{E:pushfwd-weil-cart}
we obtain the cartesian diagram
\def\mathrlap{\mathpalette\mathrlapinternal}
\def\mathrlapinternal#1#2{\rlap{$\mathsurround=0pt#1{#2}$}}
$$\xymatrix@C+5mm{%
\weilr_{S'/S}(\stX)\ar[r]^{\weilr_{S'/S}(f)}\ar[d] &\weilr_{S'/S}(\stY)\ar[r]\ar[d]
& {S\mathrlap{\smash{\;\iso\FHilb^d_{S'/S}}}}\ar[d]^-{h} \\
\HilbSt^d_{\stX/S}\ar[r]^{f_*} & \HilbSt^d_{\stY/S}\ar[r]\ar@{}[ul]|\square
& \HilbSt^d_{S'/S}\ar@{}[ul]|\square}$$
and the lemma follows.
\end{proof}

\begin{proposition}\label{P:weilr-props}
Let $S$ be a scheme and let $S'\to S$ be a finite flat morphism of finite
presentation. Let $\map{f}{\stX}{\stY}$ be a morphism of stacks over $S'$. If
$f$ has one of the properties:
\begin{enumerate}
\item locally of finite presentation,\label{PI:wp-lfp}\label{PI:wp-first}
\item formally \etale{},\label{PI:wp-etale}
\item formally unramified,\label{PI:wp-unramified}
\item formally smooth,\label{PI:wp-smooth}
\item surjective and smooth,\label{PI:wp-smooth-surj}
 \label{PI:wp-last-of-functorial}
\item a closed immersion,\label{PI:wp-closed-imm}
\item an open immersion,\label{PI:wp-open-imm}
\item an isomorphism,\label{PI:wp-iso}
\item affine;\label{PI:wp-affine}
\end{enumerate}
then so has $\weilr_{S'/S}(f)$ (also see Proposition~\pref{P:weilr-props-2} for
further properties).
\end{proposition}
\begin{proof}
Let $P$ be one of the properties and assume that $f$ has $P$. Since $P$ is
Zariski-local, we can assume that $S'\to S$ has constant rank~$d$. By
Lemma~\pref{L:property-equivalence} we can further assume that $S$ is affine,
and that $\stY=S'$.

Properties \ref{PI:wp-first}--\ref{PI:wp-last-of-functorial} of
$\weilr_{S'/S}(\stX)$ are verified using the functorial characterization of
morphisms which are locally of finite
presentation~\cite[Prop.~8.14.2]{egaIV},\ \cite[Prop.~4.15]{laumon} and the
infinitesimal criteria for formally \etale{}, unramified and smooth
maps.
Property~\ref{PI:wp-closed-imm} follows from~\cite[Lem.~9.7.9.1]{egaI_NE}
and properties~\ref{PI:wp-open-imm} and~\ref{PI:wp-iso} are obvious.
For~\ref{PI:wp-affine} it is by~\ref{PI:wp-closed-imm} enough to show that
$\weilr_{S'/S}(X)$ is represented by a scheme affine over $S$ when $X$ is the
spectrum of a polynomial ring over $\sO_{S'}$. This is straight-forward. We
refer to~\cite[\S7.6, Prop.~2, pf.\ of\ Thm.~4,
Prop.~5]{raynaud-bosch-lutkebohmert_Neron-models} for details.
%
\end{proof}

\begin{theorem}[{\cite[\S7.6, Thm.~4]{raynaud-bosch-lutkebohmert_Neron-models}}]
\label{T:repr-of-weilr-AF}
Let $S$ be an affine scheme and let $\map{g}{S'}{S}$ be a finite flat morphism
of finite
presentation. Let $X\to S'$ be a morphism of schemes. If $X$ is an
\AF{}-scheme then $\weilr_{S'/S}(X)$ is an \AF{}-scheme.
\end{theorem}
\begin{proof}
Let $X=\bigcup_\alpha U_\alpha$ be an open cover of $X$ by affines such that
every finite subset of points of $X$ lies in some $U_\alpha$. By
Proposition~\pref{P:weilr-props}, we have that $\weilr_{S'/S}(U_\alpha)$ is
affine and that $\weilr_{S'/S}(U_\alpha)\to \weilr_{S'/S}(X)$ is an open
immersion. It is then easily seen that $\coprod_\alpha
\weilr_{S'/S}(U_\alpha)\to \weilr_{S'/S}(X)$ is a Zariski covering so that
$\weilr_{S'/S}(X)$ is a scheme. Moreover, $\weilr_{S'/S}(X)$ is an
\AF{}-scheme since a finite number of points in $\weilr_{S'/S}(X)$
corresponds to a morphism $S'\times_S \coprod_{i=1}^n \Spec k_i\to X$ and
this factors through one of the $U_\alpha$'s.
\end{proof}

We say that a morphism of stacks $\map{f}{\stX}{\stY}$ is \emph{algebraic} if
for every affine scheme $T$ and morphism $T\to \stY$ the stack
$\stX\times_{\stY} T$ is an algebraic stack.

\begin{theorem}\label{T:alg-of-weilr}
Let $S'\to S$ be a finite flat morphism of finite presentation between
algebraic spaces.
Let $\map{f}{\stX}{\stY}$ be a morphism of $S'$-stacks.
\begin{enumerate}
\item If $f$ is \AF{} (e.g., representable, locally quasi-finite and separated)
then $\weilr_{S'/S}(f)$ is \AF{} and, in particular, strongly
representable.\label{TI:weilr-strongrep}
\item If $f$ is representable then so is $\weilr_{S'/S}(f)$.\label{TI:weilr-rep}
\item If $f$ is algebraic then so is $\weilr_{S'/S}(f)$.\label{TI:weilr-alg}
\end{enumerate}
In particular, if $\stX$ is an algebraic space (resp.\ an algebraic stack)
then so is $\weilr_{S'/S}(\stX)$.
\end{theorem}
\begin{proof}
By Lemma~\pref{L:property-equivalence} we can assume that $S$ is affine
and $\stY=S'$. If $f$ is as in~\ref{TI:weilr-strongrep}, then $\stX$ is an
\AF{}-scheme so that $\weilr_{S'/S}(\stX)$ is an \AF{}-scheme
according to Theorem~\pref{T:repr-of-weilr-AF}.

\ref{TI:weilr-rep} If $f$ is representable, then $\stX$ is an algebraic
space. Choose an \etale{} presentation $U\to\stX$ with $U$ an \AF{}-scheme
(e.g., a disjoint union of affine schemes). By
Theorem~\pref{T:repr-of-weilr-AF} we have that $\weilr_{S'/S}(U)$ is a
scheme. Furthermore, by Proposition~\pref{P:weilr-props}
and~\ref{TI:weilr-strongrep} we have that $\weilr_{S'/S}(U)\to
\weilr_{S'/S}(\stX)$ is \etale{}, surjective and strongly representable. Thus,
by definition $\weilr_{S'/S}(\stX)$ is an algebraic space.

\ref{TI:weilr-alg} If $f$ is algebraic, then $\stX$ is an algebraic
stack. Choose a smooth presentation $U\to\stX$ with $U$ an algebraic
space. Then $\weilr_{S'/S}(U)$ is an algebraic space by~\ref{TI:weilr-rep} and
$\weilr_{S'/S}(U)\to \weilr_{S'/S}(\stX)$ is smooth, surjective and
representable by Proposition~\pref{P:weilr-props} and~\ref{TI:weilr-rep}. Thus
$\weilr_{S'/S}(\stX)$ is an algebraic stack.
\end{proof}

We now complement Proposition~\pref{P:weilr-props} with some additional
properties.

\begin{proposition}\label{P:weilr-props-2}
Let $S$ be a scheme and let $S'\to S$ be a finite flat morphism of finite
presentation. Let $\map{f}{\stX}{\stY}$ be a morphism of algebraic stacks
over $S'$ so that $\weilr_{S'/S}(f)$ is a morphism of algebraic stacks over
$S$. If $f$ has one of the properties:
\begin{enumerate}\setcounter{enumi}{9}
\item \AF{},\label{PI:wp-AF}
\item representable,\label{PI:wp-repr}
\item locally of finite type,\label{PI:wp-lft}
\item quasi-compact,\label{PI:wp-qc}
\item quasi-affine,\label{PI:wp-qaff}
\item of finite type,\label{PI:wp-ft}
\item of finite presentation,\label{PI:wp-fp}
\item monomorphism,\label{PI:wp-mono}
\item representable and separated,\label{PI:wp-repr+sep}
\item quasi-separated,\label{PI:wp-qsep}
\item separated diagonal,\label{PI:wp-sep-diag}
\item affine diagonal,\label{PI:wp-aff-diag}
\item quasi-affine diagonal,\label{PI:wp-qaff-diag}
\item unramified diagonal (i.e., relatively Deligne--Mumford);\label{PI:wp-DM}\label{PI:wp-last}
\end{enumerate}
then so has $\weilr_{S'/S}(f)$.
\end{proposition}
\begin{proof}
As before we can assume that $S$ is an affine scheme and that $\stY=S'$.
Properties~\ref{PI:wp-AF} and~\ref{PI:wp-repr} are part of
Theorem~\pref{T:alg-of-weilr}.

\ref{PI:wp-lft} Take a smooth surjective morphism $U\to \stX$ such
that $U$ is a disjoint union of affine schemes. If $f$ is locally of finite
type, then $U\to S'$ factors through a closed immersion $U\inj W$ and a
morphism $W\to S'$ which is locally of finite presentation. Thus
by~\ref{PI:wp-closed-imm},~\ref{PI:wp-lfp} and~\ref{PI:wp-smooth-surj}, it
follows that $\weilr_{S'/S}(\stX)$ is locally of finite type.

Similarly, for property~\ref{PI:wp-qc} take a smooth surjective morphism $U\to
\stX$ with $U$ affine and the quasi-compactness of $\weilr_{S'/S}(\stX)$ follows
from~\ref{PI:wp-affine}. Property~\ref{PI:wp-qaff} is the conjunction of
properties~\ref{PI:wp-open-imm},~\ref{PI:wp-affine} and~\ref{PI:wp-qc}.

As $\weilr_{S'/S}$ preserves fiber products we have that
$\weilr_{S'/S}(\Delta_f)=\Delta_{\weilr_{S'/S}(f)}$.
Property~\ref{PI:wp-mono} [resp.\ \ref{PI:wp-repr+sep}] is equivalent to the
diagonal being an isomorphism [resp.\ a closed immersion].
Property~\ref{PI:wp-qsep} is equivalent to the quasi-compactness of
the diagonal and its diagonal. Properties~\ref{PI:wp-mono}--\ref{PI:wp-DM}
thus follow from applying the Proposition to the diagonal $\Delta_f$ and its
diagonal $\Delta_{\Delta_f}$ with the properties~\ref{PI:wp-unramified},
\ref{PI:wp-closed-imm}, \ref{PI:wp-iso}, \ref{PI:wp-affine}, \ref{PI:wp-qc},
\ref{PI:wp-qaff}.

Finally, properties~\ref{PI:wp-ft} and~\ref{PI:wp-fp} follow from
properties~\ref{PI:wp-lfp}, \ref{PI:wp-lft}, \ref{PI:wp-qc}
and~\ref{PI:wp-qsep}.
\end{proof}

\begin{remark}
%
%
If $S'\to S$ is finite and \etale{} then Proposition~\pref{P:weilr-props} also
holds for the properties ``proper'', ``flat'' and
``separated''~\cite[\S7.6, Prop.~5]{raynaud-bosch-lutkebohmert_Neron-models}.
%
\end{remark}

\begin{example}
Proposition~\pref{P:weilr-props} does not hold for the property ``proper'' nor
for the property ``finite and \etale{}''. In fact, let $S$ be arbitrary and let
$S'\to S$ be a finite flat ramified cover of degree $d$. Then
$\weilr_{S'/S}(S'\amalg S')\to S$ is \etale{} and has generic rank $2^d$ but has
lower rank over the branch locus of $S'\to S$. Thus
$\weilr_{S'/S}(S'\amalg S')\to S'$ is not proper.

Similarly, if $\stX$ is a separated algebraic stack, i.e., has proper diagonal,
then $\weilr_{S'/S}(\stX)$ need not be separated unless $S'/S$ is
\etale{}. For example, let $S'/S$ be a finite flat ramified covering of degree
$d$ as before and let $\stX=\BG(S')$ where $G$ is a finite constant group. Then
$\weilr_{S'/S}(\stX)$ is an \etale{} gerbe over $S$ with generic geometric
automorphism group $G^d$ but with automorphism group of lower rank over the
points of $S$ where $S'/S$ is ramified.
\end{example}

\begin{theorem}\label{T:pushfwd}
Let $\map{f}{\stX}{\stY}$ be a morphism of algebraic stacks.
\begin{enumerate}
\item The morphism $\map{f_*}{\HilbSt^d_X}{\HilbSt^d_Y}$ is algebraic.
\label{TI:pf-algebraic}
%
%
\item If $f$ has one of the properties~\ref{PI:wp-first}--\ref{PI:wp-last} of
Propositions~\pref{P:weilr-props} and~\pref{P:weilr-props-2}, then so has $f_*$.
\label{TI:pf-properties}
\end{enumerate}
\end{theorem}
\begin{proof}
This is the conjunction of Lemma~\pref{L:property-equivalence},
Theorem~\pref{T:alg-of-weilr} and Propositions~\pref{P:weilr-props}
and~\pref{P:weilr-props-2}.
\end{proof}

\begin{theorem}\label{T:alg-of-Hom}
Let $X\to S$ be a finite flat morphism of finite presentation and let $\stY\to
S$ be an arbitrary algebraic stack. Then the stack $\stHom_S(X,\stY)$ is
algebraic. If $\stY\to S$ has one of the properties in
Propositions~\pref{P:weilr-props} and~\pref{P:weilr-props-2}, then so has
$\stHom_S(X,\stY)\to S$.
\end{theorem}
\begin{proof}
As $\stHom_S(X,\stY)=\weilr_{X/S}(X\times_S \stY)$ this follows immediately
from Theorem~\pref{T:alg-of-weilr}.
\end{proof}

\begin{remark}\label{R:weilr-along-stacks}
The argument in~\cite[\S3.3]{olsson_Hom-stacks} shows that
Theorem~\pref{T:alg-of-Hom} remains true if $X$ is an \emph{algebraic stack}
and $X\to S$ is proper, quasi-finite and flat of finite presentation.
Indeed, any such stack admits, \etale{}-locally on $S$, a finite flat
presentation.
\end{remark}

\end{section}


\begin{section}{Algebraicity of the Hilbert functor and the Hilbert stack}
Let $\map{f}{\stU}{\stX}$ be a morphism of separated algebraic stacks. Let
$\FHilb^d_{\stU\to \stX}\subseteq \FHilb^d_{\stU}$ be the subfunctor
parameterizing families $Z\inj \stU\times_S T$ such that the composition
$Z\inj \stU\times_S T\to \stX\times_S T$ is a closed immersion. Then
$$%
\vcenter{\xymatrix{%
\FHilb^d_{\stU\to\stX}\ar[r]^{f_*}\ar[d] & \FHilb^d_{\stX/S}\ar[d] \\
\HilbSt^d_{\stU/S}\ar[r]^{f_*} & \HilbSt^d_{\stX/S}\ar@{}[ul]|\square}}$$
is cartesian and the two vertical morphisms are open immersions.

\begin{theorem}\label{T:Hilb-functor-alg}
Let $\stX/S$ be a \emph{separated} algebraic stack. Then $\FHilb^d_{\stX/S}$ is
a separated algebraic space.
\end{theorem}
\begin{proof}
We can assume that $S$ is affine. Let $\map{f}{U=\coprod_\alpha
  U_\alpha}{\stX}$ be a smooth presentation such that the $U_\alpha$'s are
affine. Then $U$ is an \AF{}-scheme and $\FHilb^d_{U/S}$ is represented by a
scheme according to Theorem~\pref{T:repr-of-Hilb-AF}. As $f$ is representable,
smooth and surjective, so is $\map{f_*}{\FHilb^d_{U\to \stX}}
{\FHilb^d_{\stX/S}}$ by Theorem~\pref{T:pushfwd}. Thus
$\FHilb^d_{\stX/S}$ is an algebraic space.
That $\FHilb^d_{\stX/S}$ is separated is Lemma~\pref{L:Hilb-separated}.
\end{proof}

As for the Hilbert functor, the algebraicity of the Hilbert stack will be an
immediate consequence of Theorem~\pref{T:pushfwd} after we have
verified that the Hilbert stack of an affine scheme is algebraic.

In the finitely presented case, the following results follow from the more
general results of~\cite[\S2.1]{lieblich_Coh_stack}. In the affine case treated
below, the proofs are a matter of elementary algebra.

\newcommand{\SYM}{\mathrm{Sym}}    

\begin{lemma}\label{L:alg-structure}
Let $B$ be an $A$-algebra and let $M$ be a locally free $A$-module of finite
rank. Then there is an $A$-algebra $Q$ which represents $B$-algebra structures
on $M$. That is, for every $A$-algebra $A'$, there is a functorial one-to-one
correspondence between $B'=B\otimes_A A'$-algebra structures on $M'=M\otimes_A
A'$ and homomorphisms $Q\to A'$. If $B$ is an $A$-algebra of finite type
(resp.\ of finite presentation), then so is $Q$.
\end{lemma}
\begin{proof}
A $B'$-algebra structure on $M'$ is given by multiplication maps
$\map{\mu}{M'\otimes_{A'} M'}{M'}$, $\map{m}{B'\otimes_{A'} M'}{M'}$, and a
unit $\map{\eta}{A'}{M'}$. Such triples of maps correspond to $A$-module
homomorphisms
$$(M\otimes_A M\otimes_A M^{\vee})
    \oplus (B\otimes_A M\otimes_A M^{\vee})
    \oplus M^{\vee}\to A'$$
and are thus represented by the symmetric algebra
$$P:=\SYM\bigl((M\otimes_A M\otimes_A M^{\vee})
      \oplus (B\otimes_A M\otimes_A M^{\vee})
      \oplus M^{\vee}\bigr).$$
That the multiplication $\mu$ is commutative, associative and compatible with
$m$ and $\eta$ can be expressed as the vanishing of the $A'$-homomorphisms
\begin{align*}
\mu  - \mu \circ \tau \;:\;
  & M'\otimes_{A'} M' \to M' \\
\mu \circ (\mu \otimes \id{M'}) - \mu \circ (\id{M'}\otimes \mu) \;:\;
  & M'\otimes_{A'} M'\otimes_{A'} M' \to M' \\
\mu \circ (m\otimes \id{M'}) - m \circ (\id{B'}\otimes \mu) \;:\;
  & B'\otimes_{A'} M'\otimes_{A'} M' \to M' \\
\mu \circ (\eta\otimes \id{M'}) - \id{M'} \;:\;
  & M' \to M'
\end{align*}
where $\map{\tau}{M'\otimes_{A'} M'}{M'\otimes_{A'} M'}$ swaps the two
factors. This vanishing is represented by a quotient $Q$ of $P$ according
to~\cite[Lem.~9.7.9.1]{egaI_NE}.

If $B$ is of finite type, then clearly so is $P$ and hence $Q$. If $B$ is of
finite presentation we use a limit argument to reduce to the noetherian case
and it follows that $Q$ is of finite presentation.
\end{proof}

\begin{theorem}\label{T:alg-of-HilbSt-aff}
Let $X$ and $S$ be affine schemes. Then $\HilbSt^d_{X/S}$ is a quasi-compact
algebraic stack with affine diagonal. If $X/S$ is of finite type (resp.\ of
finite presentation) then so is $\HilbSt^d_{X/S}$.
\end{theorem}
\begin{proof}
There is a natural morphism $\HilbSt^d_{X/S}\to \BGL_d(S)$ which maps $(Z,p,q)$
to the locally free $\sO_T$-module $p_*\sO_Z$. The stack $\BGL_d(S)/S$ is a
finitely presented algebraic stack with affine diagonal.
Lemma~\pref{L:alg-structure} shows that $\HilbSt^d_{X/S}\to \BGL_d(S)$ is
represented by affine morphisms and the theorem follows.
\end{proof}


\begin{theorem}\label{T:algebraicity-of-HilbSt}
Let $\stX/S$ be an algebraic stack. Then $\HilbSt^d_{\stX/S}$ is algebraic. If
$\stX/S$ has one of the properties: quasi-compact, quasi-separated, locally of
finite presentation, locally of finite type, separated diagonal, affine
diagonal, quasi-affine diagonal; then so has~$\HilbSt^d_{\stX/S}\to S$.
\end{theorem}
\begin{proof}
By Theorem~\pref{T:alg-of-HilbSt-aff}, we have that $\HilbSt^d_{S/S}$ is an
algebraic stack over $S$ of finite presentation and with affine diagonal. The
algebraicity of $\HilbSt^d_{\stX/S}$ thus follows from
Theorem~\pref{T:pushfwd}~\ref{TI:pf-algebraic}. The properties of
$\HilbSt^d_{\stX/S}\to S$ follow from
Theorem~\pref{T:pushfwd}~\ref{TI:pf-properties}.
\end{proof}

Note that $\HilbSt^d_{\stX/S}$ is not Deligne--Mumford nor has quasi-finite
diagonal. On the other hand, it can be seen that the open substack
$\HilbSt^{d,\unram}_{\stX/S}$ is Deligne--Mumford.


\end{section}


\begin{section}{\'{E}tale families}\label{S:etale-families}
Let $S$ be a scheme and let $X/S$ be an algebraic space. The \emph{stack of
branchvarieties}~\cite{alexeev-knutson_branch-varieties} on $X$ is
the open substack of the Hilbert stack $\HilbSt_{X/S}$ parameterizing families
$(\map{p}{Z}{T},\map{q}{Z}{X})$ such that
the geometric fibers of $p$ are reduced.
For zero-dimensional families of rank $d$, this is the open substack
$\EtSt^d_{X/S}$ of $\HilbSt^d_{X/S}$ parameterizing \etale{} families of rank
$d$. If $X/S$ is separated, it is natural to also study the
subspace~$\Et^d_{X/S}$ parameterizing \etale{} families $Z\to T$ of rank $d$
such that $Z\to X\times_S T$ is a closed immersion. Following the notation
of~\cite[Cons.~6.6]{laumon} we let $\Sec^d_{X/S}$ be the open subset of
$(X/S)^d=X\times_S \dots \times_S X$ which is the complement of the diagonals.
The symmetric group $\SG{d}$ acts by permutations on $(X/S)^d$ and this action
is free over $\Sec^d_{X/S}$. We have the following descriptions of these
stacks:

\begin{theorem}
Let $X/S$ be an arbitrary algebraic space. There is a natural isomorphism
$\EtSt^d_{X/S}\to [(X/S)^d/\SG{d}]$. If $X/S$ is separated, then the open
substack $\Et^d_{X/S}$ is identified with the algebraic space
$\Sec^d_{X/S}/\SG{d}$.
\end{theorem}
\begin{proof}
We will construct canonical morphisms in both directions. Let
$(\map{p}{Z}{T},\map{q}{Z}{X})$ be a $T$-point of $\EtSt^d_{X/S}$. The
scheme $(Z/T)^d$ is \etale{} of rank $d^d$ over $T$. The diagonals of this
scheme are open and closed, and their complement $\Sec^d_{Z/T}$ is \etale{}
of rank $d!$. This can be verified over algebraically closed points where it
is trivial. The scheme $\Sec^d_{Z/T}$ is an $\SG{d}$-torsor and comes with an
$\SG{d}$-equivariant morphism $\Sec^d_{Z/T}\to (X/S)^d$. This defines a
$T$-point of $[(X/S)^d/\SG{d}]$.

Conversely, let $W\to T$ be a $T$-point of $[(X/S)^d/\SG{d}]$, i.e., let $W/T$
be an $\SG{d}$-torsor together with an $\SG{d}$-equivariant morphism $W\to
(X/S)^d$. Let $\SG{d-1}$ be the subgroup of $\SG{d}$ acting by permuting the
first $d-1$ factors of $(X/S)^d$. This group acts freely on $W$ and the
quotient $Z=W/\SG{d-1}$ is an algebraic space, \etale{} of rank $d$ over $T$.
Moreover, the composition of $W\to (X/S)^d$ with the last projection is
$\SG{d-1}$-invariant and induces a morphism $Z\to X$. We have thus constructed
a $T$-point of $\EtSt^d_{X/S}$.

It is clear that these constructions are functorial and thus defines morphisms
$\map{F}{\EtSt^d_{X/S}}{[(X/S)^d/\SG{d}]}$ and
$\map{G}{[(X/S)^d/\SG{d}]}{\EtSt^d_{X/S}}$. It is not difficult to show that
these are (quasi-)inverses. In fact, if $(Z,p,q)$ is a $T$-point of
$\EtSt^d_{X/S}$, then we have a canonical morphism $\Sec^d_{Z/T}/\SG{d-1}\to Z$
and that this is an isomorphism can be checked over algebraically closed
points. Conversely, if $W/T$ is an $\SG{d}$-torsor, then we obtain a morphism
$W\to ((W/\SG{d-1})/T)^d$ where the $i$\textsuperscript{th} factor is the
composition of $\map{\tau_{in}}{W}{W}$ and the quotient $W\to
W/\SG{d-1}$. Again, it is easily verified that $W\to ((W/\SG{d-1})/T)^d$
induces an isomorphism of $W$ onto $\Sec^d_{(W/\SG{d-1})/T}$.
\end{proof}

\begin{remark}
The universal $\SG{d}$-torsor of $[(X/S)^d/\SG{d}]$ is $(X/S)^d$. The above
isomorphism shows that $[(X/S)^d/\SG{d-1}]=[(X/S)^{d-1}/\SG{d-1}]\times_S X$ is
the universal \etale{} rank $d$ family on $\EtSt^d_{X/S}$.
\end{remark}

\end{section}


\appendix
\begin{section}{Algebraic spaces and stacks}\label{A:stacks}
A sheaf of sets $F$ on the category of schemes $\Sch$ with the \etale{}
topology is an \emph{algebraic space} if there exists a scheme $X$ and a
morphism $X\to F$ which is represented by surjective \etale{} morphisms
of schemes~\cite[D\'ef.~5.7.1]{raynaud-gruson}, i.e., for any scheme $T$
and morphism $T\to F$, the fiber product $X\times_F T$ is a scheme and
$X\times_F T\to T$ is surjective and \etale{}.
A \emph{stack} $\stX$ is a category fibered in groupoids over $\Sch$ with the
\etale{} topology satisfying the usual sheaf
condition~\cite[D\'ef.~3.1]{laumon}.

A morphism $\map{f}{\stX}{\stY}$ of stacks is \emph{representable} if for any
scheme $T$ and morphism $T\to \stY$, the $2$-fiber product $\stX\times_\stY T$
is an algebraic space. A stack $\stX$ is \emph{algebraic} if there exists a
smooth presentation, i.e., a smooth, surjective and representable morphism
$U\to \stX$ where $U$ is a scheme (or algebraic space). A morphism
$\map{f}{\stX}{\stY}$ of stacks is \emph{quasi-separated} if the diagonal
$\Delta_{\stX/\stY}$ is quasi-compact and quasi-separated, i.e., if both
$\Delta_{\stX/\stY}$ and its diagonal are quasi-compact.

We do not require that algebraic spaces and stacks are quasi-separated nor
that the diagonal of an algebraic stack is separated. The diagonal of a (not
necessarily quasi-separated) algebraic space is represented by schemes. This
follows by effective fppf-descent of monomorphisms which are locally of finite
type. Indeed, more generally the class of locally quasi-finite and separated
morphisms is an effective class in the fppf-topology (cf.\ \cite[App.]{murre},
\cite[Exp.~X, Lem.~5.4]{sga3} or~\cite[pf.\ of 5.7.2]{raynaud-gruson}).
The diagonal of an algebraic stack $\stX$ is representable. This follows
by~\cite[pf.\ of Prop.~4.3.1]{laumon} as~\cite[Cor.~1.6.3]{laumon} generalizes
to arbitrary algebraic spaces.

In particular, if $X$ is an algebraic space (resp.\ an algebraic stack) and
$S$ and $T$ are schemes (resp.\ algebraic spaces) with morphisms $S\to X$ and
$T\to X$, then $S\times_X T$ is a scheme (resp.\ an algebraic space).
\end{section}

\bibliography{hilbert_stacks}
\bibliographystyle{dary}

\end{document}